\documentclass[11pt]{amsart}

\usepackage{fullpage}
\usepackage{amsmath}
\usepackage{amsfonts}
\usepackage{amsthm}
\usepackage{amssymb}
\usepackage{enumitem}
\usepackage{hyperref}
\usepackage{microtype}

\usepackage[backgroundcolor=lightgray]{todonotes}

\theoremstyle{plain}
\newtheorem{theorem}{Theorem}
\newtheorem{corollary}[theorem]{Corollary}
\newtheorem{proposition}[theorem]{Proposition}

\theoremstyle{definition}
\newtheorem{definition}[theorem]{Definition}
\newtheorem{example}[theorem]{Example}

\theoremstyle{remark}
\newtheorem{remark}[theorem]{Remark}

\newcommand{\beqn}{\begin{eqnarray*}}
	\newcommand{\eeqn}{\end{eqnarray*}}

\usepackage{mathtools}
\DeclarePairedDelimiter{\abs}{\lvert}{\rvert}

\DeclarePairedDelimiter{\norm}{\lVert}{\rVert}

\DeclarePairedDelimiterX{\dset}[2]{\lbrace}{\rbrace}{#1\;\delimsize|\;#2}

\newcommand{\ball}{{\overline{B}}}

\title[Hyperconvexity in partial metric spaces: challenges and outlooks]{Hyperconvexity in Partial Metric Spaces:\\ Challenges and Outlooks}

\author[D.Bugajewski]{Dariusz Bugajewski}
\address{D. Bugajewski, Department of Nonlinear Analysis and Applied Topology\\
  Faculty of Mathematics and Computer Science\\
  Adam Mickiewicz University, Pozna\'n\\
  ul.\ Umultowska 87\\
  61-614 Pozna\'n\\
  Poland}
\email{ddbb@amu.edu.pl}

\author[P.Kasprzak]{Piotr Kasprzak} 
\address{P. Kasprzak, Department of Nonlinear Analysis and Applied Topology\\
  Faculty of Mathematics and Computer Science\\
  Adam Mickiewicz University, Pozna\'n\\
  ul.\ Umultowska 87\\
  61-614 Pozna\'n\\
  Poland}
\email{piotr.kasprzak@amu.edu.pl} 

\author[O. Olela-Otafudu]{Olivier Olela-Otafudu} 
\address{O. Olela-Otafudu, Department of Mathematics and Applied Mathematics\\ University of Limpopo, Private bag X1106, Sovenga, 0727, 
South Africa}
\email{olivier.olela-otafudu@ul.ac.za}

\date{\today}

\begin{document}

\begin{abstract} 
In this article, we present several different ways to define hyperconvexity in partial metric spaces. In particular, we show that the analogue of the Aronszajn--Panitchpakdi notion of hyperconvexity fails to exhibit certain key properties present in the classical metric setting.  
\end{abstract}

\maketitle
	
\section{Introduction}

The concept of a hyperconvex metric space was introduced by Aronszajn and Panitchpakdi in~\cite{AP}, although Aronszajn had already investigated it earlier in his unpublished thesis~\cite{A}. There are several reasons why hyperconvexity is particularly interesting. First, these spaces provide a natural setting in which a metric analogue of the Hahn--Banach extension theorem holds. Second, bounded hyperconvex spaces enjoy the fixed point property for non-expansive mappings -- a result independently established by Sine~\cite{Si} and Soardi~\cite{So}, and later proven in full generality by Baillon~\cite{Ba}. Finally, from a topological point of view, every hyperconvex metric space is an absolute retract under a non-expansive retraction (see~\cite{AP}), a property that plays an important role in the fixed point theory. Further information on hyperconvexity in the metric setting can be found, for example, in \cite{Bo, BBK, EK}. 
	
Let us now turn our attention to partial metric spaces. Although these spaces are generally non-Hausdorff, as indicated in the pioneering work of Matthews~\cite{M}, many metric-like tools and concepts can be naturally extended to this setting. This can be further illustrated by the result in~\cite{BMW} which shows that in partial metric spaces compactness and sequential compactness are equivalent. Regarding their applications, partial metric spaces are most commonly associated with the study of denotational semantics of programming languages. However, they also find interesting applications in other areas, such as the geometry of normed spaces (see~\cite{ROSP2}), the domain of words, and complexity spaces (see~\cite{ROSP}). Finally, we also note that the fixed point theory in partial metric spaces is by now well-established (see, for instance,~\cite{BM, BMW, HRS, IPR} and the references therein). 

As the title of this note suggests, a question arises: how can the notion of hyperconvexity be extended to partial metric spaces? The question is all the more natural given the article~\cite{ST} of Smyth and Tsaur, in which they expressed the view that the methods they developed in their study of hyperconvexity in semi-metric spaces could be also applied to the setting of (O'Neill) partial metrics (see \cite[p.~799]{ST}). Furthermore, in the concluding remarks, they emphasized that hyperconvexity in partial metric spaces, together with related topics, calls for a more systematic investigation (see \cite[Section~7]{ST}). Motivated by these observations and remarks of Smyth and Tsaur, as well as by the fact that -- despite apparent interest within parts of the mathematical community, as reflected in conference discussions -- this problem has not yet been addressed, we attempt to answer the question posed at the beginning of this paragraph. Unfortunately, the situation turns out to be less straightforward than one might expect. No single, unambiguous definition presents itself, and multiple different approaches are possible. In particular, we show that defining hyperconvexity for partial metric spaces by direct analogy with the classical definition of Aronszajn and Panitchpakdi may lead to the loss of some of its characteristic properties.  

The paper is organized as follows. In Section 2, we recall some basic definitions and facts related to partial metric spaces that will be used throughout the paper. In Section 3, we introduce and compare several notions of hyperconvexity in the setting of partial metric spaces and examine their basic properties. Finally, Section 4 explores suitable definitions of Lipschitz continuity in partial metric spaces and their connection with a metric analogue of the Hahn--Banach extension theorem proven by Aronszajn and Panitchpakdi.

\section{Preliminaries} 
In this section, we collect some key notions and results related to partial metric spaces. Further information and examples can be found, for instance, in~\cite{BW}. 

Let $U$ be a non-empty set. A mapping $p \colon U\times U \to [0,+\infty)$ is called a \emph{partial metric} (or \emph{pmetric}) on $U$, if for all $x,y,z\in U$ the following four conditions are satisfied:
\begin{enumerate}[label=\textup{(P\arabic*)}]
 \item $x=y$ if and only if $p(x,x)=p(y,y)=p(x,y)$,
 
 \item $p(x,x)\leq p(x,y)$,

 \item $p(x,y)=p(y,x)$,

 \item $p(x,y)\leq p(x,z)+p(z,y)-p(z,z)$.
\end{enumerate}
The pair $(U,p)$ is called a \emph{partial metric space}, and $p(x,x)$ is referred to as the \emph{size} of $x$. Clearly, any metric space is a partial metric space with $p(x,x)=0$ for all $x\in U$.

Interestingly, in some sense, the converse is also true: to every partial metric space $(U,p)$ one can naturally associate a metric space, whose properties closely mirror those of $(U,p)$. For all $x,y\in U$ the corresponding metrics $p^m$ and $d_m$ are defined by
\begin{align*}
  p^m(x,y)&:=2p(x,y)-p(x,x)-p(y,y),\\
	d_m(x,y)&:=\max\bigl\{p(x,y)-p(x,x),\ p(x,y)-p(y,y)\bigr\}.
\end{align*}
For example, it is well-known that the completeness of a partial metric space $(U,p)$ -- a somewhat technical notion due to the possibility of points having positive sizes -- is equivalent to the completeness of the metric space $(U, p^m)$ or $(U, d_m)$ (see, for example,~\cite[p.~194]{M}). Yet another way to associate a metric to a partial metric was proposed by Hitzler and Seda in~\cite{HS} (see also~\cite{HRS}). They showed that the mapping $D \colon U \times U \to [0,+\infty)$ defined by 
\[
D(x,y) :=
\begin{cases}
p(x,y), & \text{if $x \neq y$,}\\
0, & \text{if $x = y$,}
\end{cases}
\]
is indeed a metric on $U$.

In a partial metric space $(U,p)$, the set $\dset{y \in U}{p(x,y)\leq r}$ may be empty for some $r>0$ and $x \in U$. To address this, we define the closed ball centered at $x\in U$ with radius $r>0$ as
\[
\ball_p(x,r):=\dset{y \in U}{p(x,y)\leq p(x,x)+r}.
\]
When $p$ is a metric, this definition clearly agrees with the usual notion of a closed ball. 

\section{Hyperconvexity}

Before we turn to the main part of our discussion, let us begin by recalling the notion of hyperconvexity introduced by Aronszajn and Panitchpakdi (see~\cite{AP}).

\begin{definition}\label{def:hyperconvex_metric_space}
A metric space $(H,d)$ is said to be \emph{hyperconvex} if for any family of closed balls $\{\ball_d(x_i,r_i)\}_{i\in I}$ in $H$ satisfying the condition $d(x_i, x_j)\leq r_i+r_j$ for any $i,j\in I$, there exists $x \in H$ with $d(x,x_i)\leq r_i$ for all $i \in I$.
\end{definition} 

When attempting to reformulate the notion of hyperconvexity from the metric to the partial metric setting, some care is required. The most straightforward -- and admittedly na\"\i{}ve -- approach (simply to replace $d$ with $p$) turns out to yield no genuinely new concept. This is due to the following simple but important observation.

\begin{proposition}\label{prop:partial_hyperconvexity_trivial}
Let $(U,p)$ be a partial metric space consisting of at least two points. Assume that for any pair of closed balls $\ball_p(x,r)$ and $\ball_p(y,R)$ whose centers satisfy the condition $p(x,y)\leq r+R$, there exists a point $z \in U$ such that $p(x,z)\leq r$ and $p(z,y) \leq R$. Then, $p$ is, in fact, a metric.
\end{proposition}

\begin{proof}
Let $x,y$ be two distinct points in $U$. Then, necessarily, $p(x,y)>0$. Now, let us fix an arbitrary parameter $\alpha \in (0,1)$, and define $r_{\alpha}:=\alpha p(x,y)$ and $R_{\alpha}:=(1-\alpha)p(x,y)$. Consider the pair of closed balls $\ball_p(x,r_{\alpha})$ and $\ball_p(y,R_{\alpha})$. Clearly, $p(x,y)=r_{\alpha}+R_{\alpha}$. Hence, by assumption, there exists a point $z_{\alpha} \in U$ such that $p(x,z_{\alpha})\leq r_{\alpha}$ and $p(z_{\alpha},y)\leq R_{\alpha}$. It follows that $ p(x,y) \leq p(x,z_{\alpha})+p(z_{\alpha},y) - p(z_{\alpha},z_{\alpha}) \leq p(x,y) - p(z_{\alpha},z_{\alpha})$, and therefore $p(z_{\alpha},z_{\alpha})=0$. In particular, this also implies that $p(x,z_{\alpha})= r_{\alpha}$ and $p(z_{\alpha},y)= R_{\alpha}$. Finally, note that if $z_\alpha=z_\beta$ for some $\alpha,\beta \in (0,1)$, then $\alpha p(x,y)=p(x,z_\alpha)=p(x,z_\beta)=\beta p(x,y)$ and, consequently, $\alpha=\beta$. 

Suppose now that $p(x,x)>0$, and choose distinct $\alpha,\beta \in (0,1)$ such that $ \alpha p(x,y) < \frac{1}{2} p(x,x)$ and $\beta p(x,y) < \frac{1}{2}p(x,x)$. Then, 
\[
 0<p(z_\alpha,z_\beta) \leq p(z_\alpha,x) + p(x,z_\beta) - p(x,x) = (\alpha+\beta)p(x,y) - p(x,x)<0,
\]
which is a contradiction. Therefore, $p(x,x)=0$. Since the point $x$ was arbitrary, we conclude that $p$ is a metric.
\end{proof} 

Hyperconvexity of a metric space means that for any family of closed balls $\{\ball_d(x_i,r_i)\}_{i\in I}$ whose centers are not too far apart, there exists a point $x$ lying in every ball $\ball_d(x_i,r_i)$ centered at $x_i$; equivalently, $\bigcap_{i\in I} \ball_d(x_i, r_i)\neq \emptyset$. There is, however, another equally valid perspective. The non-emptiness of $\bigcap_{i \in I} \ball_d(x_i,r_i)$ also means that there exists a point $y$ for which each $x_i$ lies in the ball $\ball_d(y,r_i)$. Rather than viewing $y$ as a common element of all balls centered at the $x_i$, we may instead view it as the center of balls that contain the respective points $x_i$. In metric spaces these two viewpoints coincide -- they are simply two sides of the same coin. In partial metric spaces, however, this equivalence may fail. This observation leads to the following two definitions.

\begin{definition}\label{def:hyperconvexity_ver1}
A partial metric space $(U,p)$ is said to be \emph{hyperconvex in the sense of Aronszajn and Panitchpakdi} (or, simply, \emph{AP-hyperconvex}) if for any family of closed balls $\{\ball_p(x_i,r_i)\}_{i\in I}$ in $U$ satisfying the condition $p(x_i, x_j)\leq r_i+r_j$ for all $i,j \in I$, there exists a point $x \in U$ such that
\begin{equation}\label{eq:condition_hyperconvexity_ver1}
 \text{$p(x,x_i) \leq p(x_i,x_i) + r_i$ for every $i \in I$.}
\end{equation}
\end{definition}

\begin{remark}
In Definition~\ref{def:hyperconvexity_ver1}, the radii may be taken to be non-negative rather than strictly positive. Indeed, if $r_j=0$ for some $j \in I$, then condition~\eqref{eq:condition_hyperconvexity_ver1} is satisfied by simply choosing $x:=x_j$. Since this definition directly generalizes Definition~\ref{def:hyperconvex_metric_space}, the same observation applies there as well.
\end{remark}
 
\begin{definition}\label{def:hyperconvexity_ver2}
A partial metric space $(U,p)$ is said to be \emph{nodally hyperconvex} if for any family of closed balls $\{\ball_p(x_i,r_i)\}_{i\in I}$ in $U$ satisfying the condition $p(x_i, x_j)\leq r_i+r_j$ for all $i,j \in I$, there exists a point $x \in U$ such that
\begin{equation}\label{eq:condition_hyperconvexity_ver2}
 \text{$p(x,x_i) \leq p(x,x) + r_i$ for every $i \in I$.}
\end{equation}
\end{definition}

\begin{remark}\label{rem:example_both_AP_nodal}
Note that, given a family of closed balls $\{\ball_p(x_i,r_i)\}_{i\in I}$ in a partial metric space $(U,p)$ satisfying $p(x_i,x_j)\leq r_i+r_j$ for all $i,j\in I$, the point $x$ appearing in Definitions~\ref{def:hyperconvexity_ver1} and~\ref{def:hyperconvexity_ver2} need not be the same. 

For example, let $U := \{a,b\}$ with $a$ and $b$ distinct, and equip it with the partial metric $p\colon U\times U \to [0,+\infty)$ defined by
\[
p(a,a):=2,\qquad p(a,b):=p(b,a):=2,\qquad p(b,b):=0.
\]
It is easy to verify that $(U,p)$ is both AP- and nodally hyperconvex. (This also follows from the proof of Proposition~\ref{prop:hyperconvex_extension}.) However, for the family consisting of the two balls $\ball_p(a,1)$ and $\ball_p(b,1)$, the point $x$ in Definition~\ref{def:hyperconvexity_ver1} is $b$, whereas in Definition~\ref{def:hyperconvexity_ver2} it is $a$.
\end{remark}

Somewhat surprisingly, the two proposed notions of hyperconvexity in partial metric spaces differ significantly: there exist spaces that satisfy one definition but not the other.

\begin{example}\label{ex:nodally_not_AP}
Let $U := \{a,b,c\}$ with all points distinct, and endow it with the partial metric $p\colon U\times U \to [0,+\infty)$ defined by
\[
 p(a,a):=p(a,b):=p(a,c):=3, \qquad p(b,b):=p(c,c):=0, \qquad p(b,c):=2,
\]
where the remaining values follow by symmetry. Then, for \emph{any} family of closed balls in $(U,p)$, condition~\eqref{eq:condition_hyperconvexity_ver2} is satisfied with $x := a$, implying that $(U,p)$ is nodally hyperconvex. However, if we consider the family consisting of the two balls $\ball_p(b,1)$ and $\ball_p(c,1)$, condition~\eqref{eq:condition_hyperconvexity_ver1} is not satisfied for any $x \in U$. Therefore, $(U,p)$ is not AP-hyperconvex.
\end{example}

\begin{example}\label{ex:AP_not_nodally}
Once again let us take $U := \{a,b,c\}$ with all points distinct. This time, however, let us endow it with the partial metric $p\colon U\times U \to [0,+\infty)$ defined by
\[
 p(a,a):=p(c,c):=10, \qquad p(b,b):=0, \qquad p(a,c):=30, \qquad p(a,b):=p(b,c):=15,
\]
where the remaining values follow by symmetry. 

The fact that $(U,p)$ is not nodally hyperconvex becomes evident if we consider the family of three closed balls: $\ball_p(a,19)$, $\ball_p(b,10)$, and $\ball_p(c,11)$.

We now show that $(U,p)$ is AP-hyperconvex. Let $\{\ball_p(x_i,r_i)\}_{i\in I}$ be an arbitrary family of closed balls in $U$ satisfying $p(x_i, x_j)\leq r_i+r_j$ for all $i,j \in I$. We may assume that at least one of the points $a$ or $c$ appears among the centers; otherwise, condition~\eqref{eq:condition_hyperconvexity_ver1} holds trivially  with $x:=b$. Suppose $a$ is one of the centers, and let $A := \dset{i \in I}{x_i = a}$. Then, for each $i \in A$, we obtain $10 = p(a,a)=p(x_i,x_i)\leq 2r_i$, so $r_i \geq 5$. Consequently, for any $i \in A$ we have
\[
 p(b,x_i) = p(b,a) = 15 \leq 10 + r_i = p(x_i,x_i) + r_i.
\]
Similarly, if the set $C:=\dset{i \in I}{x_i = c}$ is non-empty, we also have $p(b,x_i) \leq p(x_i,x_i) + r_i$ for all $i \in C$. From this, the AP-hyperconvexity of $(U,p)$ follows immediately.
\end{example}

\begin{remark}\label{rem:hyperconvex_total_convex}
Recall that every hyperconvex metric space $(H,d)$ is totally (Menager) convex, that is, for any points $x,y \in H$ and any non-negative numbers $\lambda,\mu$ such that $\lambda+\mu=d(x,y)$ there exists a point $z \in H$ such that $d(x,z)=\lambda$ and $d(z,y)=\mu$ (see, for example,~\cite[Theorem~2, p.~417]{AP} or~\cite[Proposition~3.1.18]{BBK}). For a partial metric space $(U,p)$, total convexity may be introduced in two natural ways. One approach is to replace $d$ with $p$ in the original definition. Another is to declare $(U,p)$ totally convex when its associated metric space $(U,p^{m}$) (or $(U,d_m)$) is totally convex in the classical sense. Regardless of which definition is adopted, Examples~\ref{ex:nodally_not_AP} and~\ref{ex:AP_not_nodally} show that partial metric spaces that are hyperconvex in either of the above senses need not be totally convex.
\end{remark}

The reasoning from the examples above extends naturally to a more general setting, providing a way to construct new hyperconvex partial metric spaces from existing ones. Before proceeding, let us recall that a partial metric space $(U,p)$ is \emph{bounded} if it is contained in some closed ball, that is, there exist an element $z \in U$ and a real number $r>0$ such that $U\subseteq \ball_p(z,r)$. As in the classical metric setting, a partial metric space $U$ is bounded if and only if $\sup_{x,y \in U} p (x,y)< +\infty$.  

\begin{proposition}\label{prop:hyperconvex_extension}
Let $(U,p)$ be a bounded partial metric space which is AP-hyperconvex, and let $a$ be a point that does not belong to $U$. Then, the set $W:=U\cup \{a\}$ admits an AP-hyperconvex partial metric $q$ that agrees with $p$ on $U\times U$.
\end{proposition}

\begin{proof}
Set $M:=1+\sup_{x,y\in U}p(x,y)$. Since $U$ is bounded, the constant $M$ if finite. Define a function $q \colon W \times W \to [0,+\infty)$ by the formula
\[
q(x,y):=\begin{cases}
p(x,y), & \text{if $x,y \in U$,}\\
M, & \text{if $x=a$ or $y=a$.}
\end{cases}
\]
It is straightforward to check that $q$ is a partial metric on $W$. 

Now, we will show that $W$ is AP-hyperconvex. Let $\{\ball_q(x_i,r_i)\}_{i \in I}$
be an arbitrary family of closed balls in $W$ satisfying the condition $q(x_i,x_j)\leq r_i+r_j$ for all $i,j \in I$. If all centers $x_i$ lie in $U$, then
\[
 \bigcap_{i \in I} \ball_p(x_i,r_i) \subseteq \bigcap_{i \in
I} \ball_q(x_i,r_i),
\]
and the former intersection is non-empty by the AP-hyperconvexity of $(U,p)$; here, $\ball_p(x_i,r_i)$ stands for the closed ball in $(U,p)$. Thus, we may assume that the point $a$ appears among the centers. We may also assume that not all centers are equal to $a$. Hence, both sets $A:=\dset {i \in I}{x_i=a}$ and $I\setminus A$ are non-empty. Observe that for each $i \in A $ and any $x \in W$ we have $
q(x_i,x)=M=q(x_i,x_i)\leq q(x_i,x_i)+r_i$. Therefore, $\ball_q(x_i,r_i)=W$ for every $i \in A$. Consequently,
\[
\bigcap_{i \in I} \ball_q(x_i,r_i) = \bigcap_{i \in I\setminus A} \ball_q
(x_i,r_i) \supseteq \bigcap_{i \in I\setminus A} \ball_p(x_i,r_i),
\]
and the last intersection is non-empty by the AP-hyperconvexity of
$(U,p)$. This completes the proof.
\end{proof}

Using the same partial metric $q$ as in the proof of Proposition~\ref{prop:hyperconvex_extension} and observing that $\ball_q(a,r)=W$ for any $r>0$, we immediately obtain the corresponding result for nodally hyperconvex partial metric spaces. Importantly, in this case there is no need to assume anything about hyperconvexity of $(U,p)$.

\begin{proposition}\label{prop:hyperconvex_extension2}
Let $(U,p)$ be a bounded partial metric space, and let $a$ be a point that does not belong to $U$. Then, the set $W:=U\cup \{a\}$ admits a nodally hyperconvex partial metric $q$ that agrees with $p$ on $U\times U$.
\end{proposition}

Beginning with a trivial (partial) metric on a singleton and applying Proposition~\ref{prop:hyperconvex_extension} or~\ref{prop:hyperconvex_extension2} repeatedly, we obtain the following corollary.

\begin{corollary}\label{cor:finite_set_hyperconvex}
Every finite set admits a partial metric that is both AP- and nodally hyperconvex.
\end{corollary}

Another fundamental property of hyperconvex metric spaces is their completeness. We now turn to examining how this property carries over to the setting of partial metric spaces. To begin, let us consider the following example.

\begin{example}\label{ex3}
Let $(E,\norm{\cdot})$ be a normed space. Following~\cite[Theorem~5.2]{M}, we can define a partial metric $p$ on $E$ by
\[
p(x,y):=\frac{\norm{x-y}+\norm{x}+\norm{y}}{2}, \quad x,y \in E.
\] 
With this definition, we have $p(x,0)=\norm{x}=p(x,x)\leq p(x,x)+r$ for any $x \in E$ and $r>0$. This shows that the partial metric space $(E,p)$ is AP-hyperconvex. 
\end{example}

\begin{remark}
One consequence of Example~\ref{ex3} is that AP-hyperconvex partial metric spaces -- unlike their classical metric counterparts -- are not necessarily complete.

The same lack of completeness can occur for nodally hyperconvex partial metric spaces. To demonstrate this, consider the open interval $U:=(0,1)$ endowed with the standard Euclidean metric $d$, and let $a:=2$. By Proposition~\ref{prop:hyperconvex_extension2} (and the proof of Proposition~\ref{prop:hyperconvex_extension}), the set $W:=(0,1)\cup\{2\}$ admits a nodally hyperconvex partial metric $q$ that coincides with $d$ on $U \times U$ and satisfies $q(2,x)=2$ for $x\in U$. It is also straightforward to check that the associated metric $q^m$ on $W$, defined by
\[
  q^m(x,y):=2q(x,y)-q(x,x)-q(y,y), \quad x,y \in W,
\]
agrees with $2d$ on $U \times U$. Since $(W,q^m)$ is not complete as a metric space, the partial metric space $(W,q)$ is likewise incomplete.
\end{remark}

\begin{remark}\label{rem:ex3_add}
Example~\ref{ex3} also shows that \emph{any} complex Banach space endowed with the partial metric defined there is AP-hyperconvex. This stands in sharp contrast to the metric case, where it is known that \emph{no} complex Banach space is hyperconvex.
\end{remark}

So far, we have considered two ways of extending hyperconvexity to the partial metric setting, both of which generally lack two key features present in the classical metric case: total convexity and completeness. (The third feature, related to the fixed-point property for non-expansive mappings, will be discussed later in Section~\ref{sec:lipschitz}). As is commonly done, a partial metric space is said to possess a property if one of its associated metric spaces has that property. We will adopt this approach here, leading to the following definition.

\begin{definition}\label{hpms}
A partial metric space $(U,p)$ is called \emph{$p^m$-hyperconvex} \textup(respectively, \emph{$d_m$-hyperconvex}, or \emph{$D$-hyperconvex}\textup) when the associated metric space $(U,p^m)$ \textup(respectively, $(U,d_m)$, or $(U,D)$\textup) is hyperconvex in the classical sense. 
\end{definition} 

Now, we are going to establish the connections between the notions of hyperconvexity introduced above and those given in Definitions~\ref{def:hyperconvexity_ver1} and~\ref{def:hyperconvexity_ver2}. 

\begin{remark}
A similar argument to the one presented in the proof of Proposition~\ref{prop:partial_hyperconvexity_trivial} shows that if $(U,p)$ is a partial metric space with at least two points which is additionally $D$-hyperconvex, then $p$ is a metric. Consequently, all the notions of hyperconvexity in partial metric spaces considered in this paper coincide in this case.
\end{remark}

\begin{proposition}\label{prop:pm_implies_p}
Let $(U,p)$ be a partial metric space. If $U$ is $p^m$-hyperconvex, then it is also nodally hyperconvex.
\end{proposition}

\begin{proof}
Let $\{\ball_p(x_i,r_i)\}_{i \in I}$ be an arbitrary collection of closed balls in $(U,p)$ such that
\[
 \text{$p(x_i,x_j)\leq r_i+r_j$ for all $i,j \in I$.}
\]
In particular, we have $p(x_i,x_i)\leq 2r_i$ for $i \in I$. Then, for any $i,j \in I$ we clearly get
\[
 p^m(x_i,x_j)=2p(x_i,x_j)-p(x_i,x_i)-p(x_j,x_j) \leq \bigl[2r_i - p(x_i,x_i)\bigr] + \bigl[2r_j - p(x_j,x_j)\bigr]. 
\]
Hence, by $p^m$-hyperconvexity of $U$, there exists a point $x \in U$ such that
\[
 \text{$p^m(x,x_i) \leq 2r_i-p(x_i,x_i)$ for all $i \in I$,} 
\]
which in turn implies that
\[
 \text{$p(x,x_i)\leq \frac{1}{2}p(x,x)+r_i$ for all $i \in I$.}
\]
This proves that $U$ is nodally hyperconvex.
\end{proof}

\begin{proposition}\label{prop:dm_implies_p}
Let $(U,p)$ be a partial metric space. If $U$ is $d_m$-hyperconvex, then it is both AP- and nodally hyperconvex.
\end{proposition}

\begin{proof}
Let $\{\ball_p(x_i,r_i)\}_{i \in I}$ be an arbitrary collection of closed balls in $(U,p)$ such that
\[
 \text{$p(x_i,x_j)\leq r_i+r_j$ for all $i,j \in I$.}
\]
Then, for any $i,j \in I$ we get
\[
p(x_i,x_j)-p(x_i,x_i) \leq p(x_i,x_j) \leq r_i + r_j.
\]
And, similarly, 
\[
p(x_i,x_j)-p(x_j,x_j) \leq r_i + r_j,
\]
which shows that
\[
 \text{$d_m(x_i,x_j) \leq r_i + r_j$ for all $i,j \in I$.}
\]
Hence, by $d_m$-hyperconvexity of $U$, there exists a point $x \in U$ such that
\[
 \text{$d_m(x,x_i) \leq r_i$ for all $i \in I$.} 
\]
Therefore,
\[
 \text{$p(x,x_i) \leq p(x,x)+r_i$ for all $i \in I$} 
\]
and
\[
 \text{$p(x,x_i) \leq p(x_i,x_i)+r_i$ for all $i \in I$.} 
\]
This completes the proof.
\end{proof}

\begin{remark}\label{rem2}
Note that Propositions~\ref{prop:pm_implies_p} and~\ref{prop:dm_implies_p} are not reversible. Indeed, the partial metric space $(U,p)$ considered in Remark~\ref{rem:example_both_AP_nodal} is both AP- and nodally hyperconvex, but is neither $d_m$-, or $p^m$-hyperconvex, because hyperconvex metric spaces are totally convex and, as a consequence, cannot be finite unless they consist of a single point (see Remark~\ref{rem:hyperconvex_total_convex}).
\end{remark}

\begin{remark}\label{rem:normed_space_AP}
In connection with the previous remark, we note that there also exist AP-hyperconvex partial metric spaces that are not $p^m$-hyperconvex, even though their associated metric spaces are complete and totally convex. For example, consider the Banach space of all null real sequences $U:=c_0$ with the supremum norm $\norm{\cdot}_\infty$. As observed in Example~\ref{ex3}, $c_0$ with the partial metric
 \[
p(x,y):=\frac{\norm{x-y}_\infty+\norm{x}_\infty+\norm{y}_\infty}{2}, \quad x,y \in c_0,
\] 
is AP-hyperconvex. Its associated metric $p^m(x,y)=\norm{x-y}_\infty$ coincides with the standard metric on $c_0$, which is known not to be hyperconvex (see, e.g., \cite[Example~3.1.11]{BBK}), even though it is complete and totally convex. 
\end{remark}

The discussion so far indicates that a partial metric space that is $p^m$- or $d_m$-hyperconvex (or both) need not be $D$-hyperconvex. This fact can be illustrated directly with the following simple example.

\begin{example}\label{ex1}
Let $U:=\mathbb R$ be endowed with the partial metric $p$ defined by $p(x,y)=1+\abs{x-y}$ for $x,y \in \mathbb R$. Its associated metrics are easily seen to be $p^m(x,y)=2\abs{x-y}$ and $d_m(x,y)=\abs{x-y}$ for $x,y\in \mathbb R$. Therefore, $(\mathbb R,p)$ is both $p^m$- and $d_m$-hyperconvex, since $\mathbb R$ with the standard Euclidean metric (or any of its positive multiples) is hyperconvex.
 
Clearly, $(\mathbb R,p)$ is not $D$-hyperconvex, since $p$ is not a metric. To see this directly, it suffices to take two closed balls: $\ball_D(1,1)=\{1\}$ and $\ball_D(2,1)=\{2\}$. 
\end{example} 

More surprisingly, a partial metric that is $p^m$-hyperconvex need not be $d_m$-hyperconvex.

\begin{example}\label{ex:pm_but_not_dm}
It is well-known that the space $\mathbb R^3$, endowed with the metric $d$ induced by the $l^1$-norm $\norm{\cdot}_1$, is not hyperconvex (see~\cite[Theorem~3.2.2]{BBK} or~\cite{BBP}). Interestingly, this changes when the metric $d$ is restricted to the subset
\[
 H:=\dset[\big]{(t,0,0) \in \mathbb R^3}{t \in [0,1]} \cup \dset[\big]{(0,t,0) \in \mathbb R^3}{t \in [0,1]} \cup \dset[\big]{(0,0,t) \in \mathbb R^3}{t \in [0,1]}
\]
(see~\cite[Example~3.4.17]{BBK}). On this subset, we define a partial metric $p$ by
\[
p(x,y):=\frac{\norm{x-y}_1+\norm{x}_1+\norm{y}_1}{2}, \quad x,y \in H, 
\] 
which turns $(H,p)$ into a $p^m$-hyperconvex space, since the associated metric $p^m$ coincides with $d$. However, $(H,p)$ fails to be $d_m$-hyperconvex. To see why, consider the points $e_1:=(1,0,0)$, $e_2:=(0,1,0)$ and $e_3:=(0,0,1)$, and the corresponding closed balls $\ball_{d_m}(e_i,\frac{1}{2})$. A straightforward computation shows that $d_m(e_i,e_j) = 1$ for all distinct $i,j \in \{1,2,3\}$. Now, suppose there exists a point $x\in H$ that lies in all three balls. Then, we would have $d_m(x,e_i)\leq \frac{1}{2}$ for each $i=1,2,3$. By symmetry, we may assume that $x:=(s,0,0)$ for some $s \in [0,1]$. But this immediately gives $d_m(x,e_2)=1$. This contradiction shows that no such point exists, confirming that $(H,p$) is not $d_m$-hyperconvex.
\end{example}

\section{Lipschitz continuity}
\label{sec:lipschitz}

The question of how to properly define a Lipschitz continuous mapping in a partial metric space $(U,p)$ is a subtle one. It is closely related to the notion of hyperconvexity, as Aronszajn and Panitchpakdi showed that any Lipschitz continuous map $f \colon A \to H$, defined on a non-empty subset $A$ of a metric space $X$ with values in a hyperconvex metric space $H$, can be extended to the whole space $X$ while preserving its Lipschitz constant (cf.~\cite[Theorem~3]{AP}); interestingly, the converse also holds (cf.~\cite[Theorem~2]{AP}). In this context, it is also worth mentioning the famous Baillon's theorem, which states that bounded hyperconvex metric spaces have the fixed point property for non-expansive maps (see~\cite{Ba}).

In the paper~\cite{M}, Matthews took a natural approach. He extended the classical notion of a contraction by simply replacing the metric with a partial metric (see~\cite[Theorem~5.3]{M}). Following this idea, we say that a mapping $f \colon U \to U$ is \emph{Lipschitz continuous \textup(in the sense of Matthews\textup)} if there exists a constant $L > 0$ such that
\[
p(f(x),f(y)) \leq L p(x,y)
\]
for all $x, y \in U$. However, it is worth noting that under this definition, not all constant mappings are Lipschitz continuous -- a rather surprising fact. In some cases, it may even happen that ,,most'' constant maps fail to be Lipschitz continuous.

\begin{example}
For any $x,y \in \mathbb R$ set
\[
 p(x,y) := \frac{\abs{x-y}+\abs{x}+\abs{y}}{2}.
\]
Obviously $p$ defines a partial metric on $\mathbb R$ that is AP-, nodally and $p^m$-hyperconvex (see Example \ref{ex3} and Proposition~\ref{prop:pm_implies_p}). Now, for a fixed $c \neq 0$ let $f_c \colon \mathbb R \to \mathbb R$ be the constant function given by $f_c(x)=c$ for $x \in \mathbb R$. If $f_c$ were Lipschitz continuous in the sense of Matthews with some constant $L>0$, then we would have   
\[
 \abs{c}=p(f_c(0),f_c(0)) \leq Lp(0,0)=0,
\]  
which is clearly impossible.
\end{example}

Another definition of Lipschitz continuity in partial metric spaces can be given following the considerations in~\cite{IPR}. A mapping $f \colon U \to U$ is said to be \emph{Lipschitz continuous \textup(in the sense of Ili\'{c}--Pavlovi\'{c}--Rako\v{c}evi\'{c}\textup)} if there exists a constant $L > 0$ such that
\[
p(f(x),f(y)) \leq \max\bigl\{L p(x,y), p(x,x), p(y,y)\bigr\}
\]
for all $x, y \in U$. This definition also appears to be unsatisfactory. As shown in~\cite[Theorem~5.8]{BMW}, a constant mapping $f_c \colon U \to U$, taking the single value $c$, is Lipschitz continuous in the above sense for any $L \in [0,1)$ if and only if $c$ belongs to the so-called \emph{bottom set}
\[
 U_p:=\dset[\big]{x \in U}{p(x,x)=\inf_{y \in U}p(y,y)}.
\]
(Although $(U,p)$ is assumed to be complete in the statement of~\cite[Theorem~5.8]{BMW}, completeness is in fact not required for this result.) In particular, if the bottom set is empty, every constant mapping fails to be a contraction in the sense of Ili\'{c}--Pavlovi\'{c}--Rako\v{c}evi\'{c}.      

All of the above discussion leads us to propose yet another definition of a Lipschitz mapping in partial metric spaces. We say that a mapping  $f \colon U \to U$ is \emph{Lipschitz continuous} if there exists a constant $L>0$ such that
\begin{equation}\label{eq:L1}
p(f(x),f(y)) \leq \max\bigl\{L p(x,y), p(x,x), p(y,y), p(f(x),f(x)), p(f(y),f(y))\bigr\}
\end{equation}
for all $x,y \in U$. Note first that when $p$ is an ordinary metric, condition~\eqref{eq:L1} reduces to the classical Lipschitz condition. Furthermore, every constant mapping $f_c \colon U \to U$ is Lipschitz continuous (in the above sense) for any constant $L>0$. Indeed, for all $x,y \in U$ we have
\begin{align*}
 p(f_c(x),f_c(y))&=p(c,c) =p(f_c(x),f_c(x))\\
 &\leq  \max\bigl\{L p(x,y), p(x,x), p(y,y), p(f_c(x),f_c(x)), p(f_c(y),f_c(y))\bigr\}.
\end{align*}
Finally, observe that if $f$ is the identity mapping, then inequality~\eqref{eq:L1} holds with $L=1$, meaning that $f$ is non-expansive, just as in the classical metric setting. (In fact, it is straightforward to verify that for the identity mapping~\eqref{eq:L1} becomes an equality when $L=1$.) 

To conclude, we remark that, instead of condition~\eqref{eq:L1}, one may also consider the alternative condition
\begin{equation*}\label{eq:L2}
p(f(x),f(y)) \leq \max\bigl\{L p(x,y), p(x,x), p(y,y), p(x,f(x)), p(y,f(y))\bigr\},
\end{equation*}
or possibly a suitable combination of the two.

At the end of this section, we present an example demonstrating that in the case of bounded AP- and nodally hyperconvex partial metric spaces non-expansive mappings do not necessarily possess the fixed point property. Thus, another key feature of hyperconvex metric spaces can be lost when moving to the partial metric setting.

\begin{example}
Let $a$ and $b$ be two distinct points. It is easy to check that the set $U:=\{a,b\}$ endowed with the partial metric $p \colon U \times U \to [0,+\infty)$ given by
\[
p(a,a):=2,\quad p(a,b):=p(b,a):=3,\quad p(b,b):=2.
\]
is both AP- and nodally hyperconvex. Also, define a mapping $f \colon U \to U$ setting $f(a):=b$ and $f(b):=a$. Clearly, $f$ is fixed-point free. However,
\begin{align*}
 p(f(a),f(a))&=p(b,b)=2=p(a,a),\\
 p(f(a),f(b))&=p(b,a)=p(a,b),\\
 p(f(b),f(b))&=p(a,a)=2=p(b,b).
\end{align*}
Thus, $f$ is an isometry in the sense of Matthews (defined analogously to a Lipschitz map), and, consequently, it is an isometry in the sense Ili\'{c}--Pavlovi\'{c}--Rako\v{c}evi\'{c} as well as in the sense of~\eqref{eq:L1}. This demonstrates that, in general, a Baillon-type fixed point theorem does not hold for AP- and nodally hyperconvex partial metric spaces.
\end{example}

\section{Conclusions}
This article aimed to explore different approaches to hyperconvexity in partial metric spaces. The situation proved more subtle than in the classical metric setting, with several \emph{non-equivalent} notions of hyperconvexity. While certain connections between these notions have been established, several questions remain open. We summarize them below:
\begin{enumerate}
 \item In Example~\ref{ex3}, we showed that every normed space $E$ is AP-hyperconvex under the partial metric defined there. Does the same hold for nodal hyperconvexity?

\item Proposition~\ref{prop:pm_implies_p} shows that a partial metric space $(U,p$) that is $p^m$-hyperconvex is also nodally hyperconvex. Can such a space be AP-hyperconvex as well?

\item Example~\ref{ex:pm_but_not_dm} demonstrates that $p^m$-hyperconvexity does not necessarily imply $d_m$-hyper\-con\-ve\-xi\-ty. Does the converse hold, and if so, under what conditions?
\end{enumerate}

\end{document}